\documentclass{amsart}

\usepackage{tikz}

\title{The curve complex and covers via hyperbolic 3-manifolds}

\author{Robert Tang}
\address{Mathematics Institute, University of Warwick, Coventry, CV4 7AL, UK}
\email{robert.tang@warwick.ac.uk}
\urladdr{http://www.warwick.ac.uk/~mariam/}
\date{13 April 2011}
\thanks{}

\numberwithin{equation}{section}

\newtheorem{Thm}[equation]{Theorem}

\newtheorem{Lem}[equation]{Lemma}

\theoremstyle{definition}

\theoremstyle{remark}
\newtheorem{Rem}[equation]{Remark}

\newcommand{\R}{\mathbb{R}}
\newcommand{\Hy}{\mathbb{H}}

\newcommand{\C}{\mathcal{C}}

\newcommand{\lab}[1]{\label{#1}}

\begin{document}

	\begin{abstract}
	
	Rafi and Schleimer recently proved that the natural relation between curve complexes induced by a covering map between two surfaces is a quasi-isometric embedding. We offer another proof of this result using a distance estimate via hyperbolic 3-manifolds.
	\end{abstract}
	
		\maketitle
	
	\section{Introduction}
	
	Let $S$ be a compact orientable surface of genus $g$ with $n$ boundary components. Define a simplicial complex $\C(S)$, called the \emph{curve complex}, as follows. Let the vertex set be the free homotopy classes of essential non-peripheral simple closed curves. We have a simplex for every set of homotopy classes which can be realised simultaneously disjointly. In particular, two classes are adjacent in the {1-skeleton} of $\C(S)$ if and only if they have disjoint representatives on $S$.

	We will assume that the \emph{complexity} $\xi(S) = 3g + n -3$ is at least 2 in this paper. This will guarantee that $\C(S)$ is non-empty, connected and that all surfaces under consideration are of hyperbolic type. (There are modified definitions for the curve complex in the low-complexity cases, however, we shall not deal with them here.)

	The distance $d_S(a,b)$ between vertices $a$ and $b$ in $\C(S)$ is defined to be the length of the shortest edge-path connecting them. This can be thought of as measuring how ``complicated'' their intersection pattern is. Endowed with this distance, the curve complex has infinite diameter and is also Gromov hyperbolic \cite{MM1,bhb-int}.
	
	The curve complex has had many applications to both the study of mapping class groups and Teichm\"uller theory. More recently, it was a key ingredient of the proof of the Ending Lamination Theorem \cite{BCM-ELC, bhb-elt}.

	We will be focussed on a map between curve complexes induced by a covering map.
	
	Suppose $P: \Sigma \rightarrow S$ is finite-degree covering map between two surfaces. The preimage $P^{-1}(a)$ of a curve $a\in\C(S)$ is a disjoint union of simple closed curves on $\Sigma$. This induces a one-to-many map $\Pi:\C(S) \rightarrow \C(\Sigma)$ where $\Pi(a)$ is defined to be the set of homotopy classes of the curves in $P^{-1}(a)$.
	
	To establish some notation, given non-negative real numbers $A,B$ and $K$, write $A \prec_K B$ to mean $A \leq KB + K$. We shall write $A \asymp_K B$ if both $A \prec_K B$ and $B \prec_K A$ hold.
	
	Let $(X,d)$ and $(X',d')$ be pseudo-metric spaces and $f:X\rightarrow X'$ a (possibly one-to-many) function. We call $f$ a \emph{$K$-quasi-isometric embedding} if there exists a positive constant $K$ such that for all $x,y\in X$ we have $d(x,y) \asymp_K d'(x',y')$ whenever $x'\in f(x)$ and $y'\in f(y)$. In addition, if the $K$-neighbourhood of $f(X)$ equals $X'$ then we call $f$ a \emph{$K$-quasi-isometry} and say $X$ and $X'$ are \emph{$K$-quasi-isometric}. Furthermore, we call $K$ the \emph{quasi-isometry constant}.

	In this paper we give a new proof of the following theorem, originally due to Rafi and Schleimer \cite{saul-covers}.

	\begin{Thm}\lab{theorem}
	Let $P: \Sigma \rightarrow S$ be a covering map of degree $\deg P < \infty$. Then the map $\Pi:\C(S) \rightarrow \C(\Sigma)$ defined above is a $K$-quasi-isometric embedding, where $K$ depends only on $\xi(\Sigma)$ and $\deg P$.
	\end{Thm}

	This result can be interpreted as saying that one cannot tangle or untangle two curves on $S$ by too much when passing to finite index covers.

	Theorem \ref{theorem} was first proved in \cite{saul-covers} using Teichm\"uller theory and subsurface projections. Our approach uses an estimate for distance in the curve complex via a suitable hyperbolic 3-manifold homeomorphic to $S\times\R$ with a modified metric. This allows us to naturally compare distances by taking a covering map between the respective 3-manifolds. The estimate (Theorem \ref{estim}) arises from work towards the Ending Lamination Theorem \cite{BCM-ELC} and is made explicit in \cite{bhb-elt}. More details will be given in Section \ref{3mfld}.
	
	\begin{Rem}
	A consequence of Theorem \ref{theorem} and Gromov hyperbolicity of $\C(\Sigma)$ is that the image $\Pi(\C(S))$ is quasi-convex.
	\end{Rem}
	
	\section{The induced map between curve complexes}
	
	Let $P: \Sigma \rightarrow S$ be a finite index covering of surfaces and assume $\xi(S) \geq 2$. The cover $P$ induces a one-to-many map $\Pi: \C(S) \rightarrow \C(\Sigma)$ defined by declaring $\Pi(a)$ to be the set of homotopy classes in $P^{-1}(a)$. Observe that $P^{-1}(a)$ is a nonempty disjoint union of circles and that any component of $P^{-1}(a)$ must be essential and non-peripheral. Thus $\Pi$ sends vertices of $\C(S)$ to (non-empty) simplices of $\C(\Sigma)$.
	
	The preimages of disjoint curves in $S$ must themselves be disjoint. Therefore, if we choose a preferred vertex in each $\Pi(a)$ then we can send any edge-path in $\C(S)$ to an edge-path in $\C(\Sigma)$. (We can extend this to higher dimensional simplices if so desired.) This gives us the following distance bound.
		
	\begin{Lem}\lab{simp}
	Let $a$ and $b$ be curves in $\C(S)$ and suppose $\alpha\in\Pi(a)$ and $\beta\in\Pi(b)$. Then
	\[ d_\Sigma(\alpha,\beta) \leq d_S(a,b).\]\qedhere
	\end{Lem}

	\section{Estimates from 3-Manifolds}\lab{3mfld}
		
	We now introduce the background required to state an estimate for distances in the curve complex using a suitable hyperbolic 3-manifold.
	
	Let $S$ be a surface with $\xi(S) \geq 2$ and fix distinct curves $a$ and $b$ in $\C(S)$. By simultaneous uniformisation, there exists a hyperbolic 3-manifold ${(X,d) \cong \mathrm{int}(S) \times \R}$ with a preferred homotopy equivalence $f$ to $S$ such that the unique geodesic representatives of the two curves, denoted $a^*$ and $b^*$, are arbitrarily short. In fact, such a 3-manifold can be chosen so that there are no accidental cusps. For a reference, see \cite{bhb-tight}.
		
	Recall that the \emph{injectivity radius} at a point $x\in X$ is equal to half the infimum of the lengths of all nontrivial loops in $X$ passing through $x$. The \emph{$\epsilon$-thin part} (or just thin part) of $X$ is the set of all points whose injectivity radius is less than $\epsilon$. By the well known Margulis lemma, the thin part comprises of Margulis tubes and cusps whenever $\epsilon$ is less than the Margulis constant. The \emph{thick part} of $X$ is the complement of the thin part.
	We shall fix such an $\epsilon$ for the rest of this paper.

	Let $\Psi_\epsilon(X)$ denote the \emph{non-cuspidal part} of $X$, that is, $X$ with its $\epsilon$-cusps removed. We will write $\Psi(X)$ for brevity. Since $X$ has no accidental cusps we see that if $S$ is closed then $\Psi(X) = X$.
		
	Define the \emph{electrified length} of a path in $\Psi(X)$ to be its total length occurring outside the Margulis tubes of $X$. More formally, we take the one-dimensional Hausdorff measure of its intersection with the thick part of $X$. This induces a reduced pseudometric $\rho_X$ on $\Psi(X)$ obtained by taking the infimum of the electrified lengths of all paths connecting two given points. One can show that the infimum is attained, for example, by taking a path which connects Margulis tubes by shortest geodesic segments.

	The distance $\rho_X$ shall be referred to as the \emph{electrified distance} on $\Psi(X)$ with respect to its Margulis tubes or the reduced pseudometric obtained by electrifying the Margulis tubes.

	\begin{Thm}\cite{bhb-elt}\lab{estim}
		Let $a$ and $b$ be curves in $\C(S)$ whose geodesic representatives $a^*$ and $b^*$ in $X$ have $d$-length at most $L \geq 0$. Then \[d_S(a,b) \asymp_K \rho_X(a^*,b^*)\] where the constant $K$ depends only on $\xi(S)$, $L$ and $\epsilon$.
	\end{Thm}

	This estimate follows from the construction of geometric models for hyperbolic 3-manifolds used in the proof of the Ending Lamination Theorem. We refer the reader to \cite{bhb-elt} and \cite{BCM-ELC} for an in-depth discussion.
	
	\begin{Rem}\lab{avoid}
	For Theorem \ref{estim} to make sense we need to ensure that $a^*$ and $b^*$ are indeed contained in $\Psi(X)$. This can be done using a pleated surfaces argument, such as in \cite{bhb-length}, to show that closed geodesics of bounded length in $X$ avoid cusps provided $\epsilon$ is sufficiently small.
	\end{Rem}
		
	\section{Proof of Theorem \ref{theorem}}
	
	\subsection{Closed surfaces}\lab{closed}
	
	We first prove the main theorem for closed surfaces. The required modifications for the general case shall be dealt with in Section \ref{boundary}.
	
	Fix a length bound $L$ and a sufficiently small value of $\epsilon$. Let $P: \Sigma \rightarrow S$ be a covering map. Fix curves $a,b$ in $\C(S)$ and choose $\alpha\in\Pi(a)$ and $\beta\in\Pi(b)$. From Lemma \ref{simp}, we have $d_\Sigma(\alpha,\beta) \leq d_S(a,b)$ so it remains to prove the reverse inequality.
	
	Let $(X,d) \cong \mathrm{int}(S) \times \R$ be a hyperbolic 3-manifold with a homotopy equivalence $f$ to $S$ as described in Section \ref{3mfld}. We also assume that $a^*$ and $b^*$ have length at most $\frac{L}{\deg P}$ in $X$. There exists a covering map ${Q: \Xi \rightarrow X}$, where ${\Xi \cong \mathrm{int}(\Sigma)\times\R}$, and a homotopy equivalence $\tilde{f}: \Xi \rightarrow \Sigma$ such that $P\circ\tilde{f} = f\circ Q$. Note that $Q(\alpha^*)=a^*$ and $Q(\beta^*)=b^*$, hence $\alpha^*$ and $\beta^*$ have length bounded above by $L$.

	Let $\rho_X$ and $\rho_\Xi$ be the respective pseudometrics on $X$ and $\Xi$ obtained by electrifying their $\epsilon$-tubes.
	
		\begin{Lem}\lab{lips}
		The map $Q$ is 1-Lipschitz with respect to $\rho_\Xi$ and $\rho_X$.
		\end{Lem}

		\proof
		Using the definition of injectivity radius and the $\pi_1$-injectivity of covering maps, we see that the $\epsilon$-thin part of $\Xi$ is sent into that of $X$. It follows that the electrified lengths of paths cannot increase under $Q$.
		\endproof
		
		This result, together with Theorem \ref{estim}, proves Theorem \ref{theorem} for closed surfaces.
				
	\subsection{Surfaces with boundary}\lab{boundary}
	
	We now assume $S$ has non-empty boundary. Recall that $\Psi_\epsilon(X)$ denotes $X$ with its $\epsilon$-cusps removed.
	
	\begin{Lem}\lab{flex}
	Let $X$ be a hyperbolic 3-manifold. Choose small constants $\delta > \delta'$ and let $\rho$ and $\rho'$ be the pseudometrics on $\Psi_\delta(X)$ and $\Psi_{\delta'}(X)$ obtained by electrifying along their respective $\epsilon$-tubes. Then the natural retraction 
\[r: (\Psi_{\delta'}(X),\rho') \rightarrow (\Psi_{\delta}(X),\rho)\] is $R$-Lipschitz, where $R = \frac{\sinh\delta}{\sinh\delta'}$.
	\end{Lem}
	
		\proof
	Begin with a geodesic arc $\gamma$ in $\Psi_{\delta'}(X)$. We will show that the length of $\gamma$ can only increase by a bounded multiplicative factor under the retraction $r$.
	
	If $\gamma$ is contained in $\Psi_\delta(X)$ then we are done. So suppose $\gamma$ meets some $\delta$-cusp $C$ of $X$. This cusp contains a $\delta'$-cusp $C'$ which does not meet $\gamma$. The cusps $C$ and $C'$ lift to nested horoballs $H$ and $H'$ in the universal cover $\Hy^3$. We can arrange for the horospheres $\partial H$ and $\partial H'$ to be horizontal planes at heights 1 and $R>1$ respectively in the upper half-space model. Using basic hyperbolic geometry and the definition of injectivity radius, we can show that $R = \frac{\sinh\delta}{\sinh\delta'}$.
	
	Now define $\pi$ and $\pi'$ to be the nearest point projections from $\Hy^3$ to $H$ and $H'$ respectively. Taking a lift $\tilde{\gamma}$ of $\gamma$, we see that
	\[\mathrm{length}(\tilde{\gamma}\cap H) \geq \mathrm{length}(\pi'(\tilde{\gamma}\cap H)) = \frac{1}{R}\mathrm{length}(\pi(\tilde{\gamma}\cap H)).\]
	Thus, by projecting each arc of $\gamma$ inside a $\delta$-cusp to the boundary of that cusp, we create a new path whose length is at most $R\times \mathrm{length}(\gamma)$.
	\endproof

	Let $Q: \Xi \rightarrow X$ be the covering map as described in Section \ref{closed}. Observe that \[\Psi_{\epsilon}(X) \subseteq Q(\Psi_\epsilon(\Xi)) \subseteq \Psi_{\epsilon'}(X)\] where $\epsilon' = \frac{\epsilon}{\deg P}$. As before, let $\rho_\Xi$ and $\rho_X$ be the pseudometrics on $\Psi_\epsilon(\Xi)$ and $\Psi_{\epsilon}(X)$ obtained by electrifying along their respective $\epsilon$-tubes.  Combining Lemma \ref{flex} with the proof of Lemma \ref{lips} gives us the following.
	
	\begin{Lem}
	Let $r: \Psi_{\epsilon'}(X) \rightarrow \Psi_{\epsilon}(X)$ be the natural retraction. Then the composition $r \circ Q: \Psi_\epsilon(\Xi) \rightarrow \Psi_{\epsilon}(X)$ is $R$-Lipschitz with respect to $\rho_\Xi$ and $\rho_X$. Moreover, the constant $R$ depends only on $\deg P$. \qed
	\end{Lem}

	Finally, ${r \circ Q(\alpha^*) = a^*}$ and ${r \circ Q(\beta^*) = b^*}$ and so the rest of the argument follows as in Section \ref{closed}.

	\section*{Acknowledgements}

	I am grateful to both Brian Bowditch and Saul Schleimer for many helpful discussions and suggestions. This work was supported by the Warwick Postgraduate Research Scholarship.
	
	\bibliography{mybib}		
	\bibliographystyle{amsplain}

\end{document}